\documentclass[11pt]{article}
\usepackage[margin=0.9in]{geometry}
\usepackage{float}
\usepackage{hyperref}
\usepackage{titling}
\usepackage{blindtext}
\usepackage{amsmath}
\usepackage{amssymb}
\usepackage{mathtools}
\usepackage[numbers]{natbib}
\usepackage{tikz}
\usepackage{graphicx}
\usepackage{siunitx}
\usepackage{gensymb}
\usepackage{url}
\usepackage{mathrsfs}
\usepackage{amsthm}
\usepackage{cleveref}
\usepackage{tkz-euclide}
\usepackage{titlesec}
\usepackage{url, verbatim}
\usepackage{enumerate}
\usepackage{pgfplots}
\usetikzlibrary{calc,patterns,quotes,shapes,arrows,through,intersections,decorations.pathreplacing}
\allowdisplaybreaks

\numberwithin{equation}{section}

\newcommand{\ubox}{\overline{\text{dim}}_\text{B}}
\newcommand{\lbox}{\underline{\text{dim}}_\text{B}}
\newcommand{\boxd}{\text{dim}_\text{B}}
\newcommand{\aso}{{\text{dim}}_\text{A}}

\newcommand{\low}{{\text{dim}}_\text{L}}
\newcommand{\haus}{{\text{dim}}_\text{H}}

\newcommand{\asospec}{{\text{dim}}^{\theta}_\text{A}}
\newcommand{\lowspec}{{\text{dim}}^{\theta}_\text{L}}

\newcommand{\lset}{L(\Gamma)}
\newcommand{\kmin}{k_{\min}}
\newcommand{\kmax}{k_{\max}}
\newcommand{\pmin}{p_{\min}}
\newcommand{\pmax}{p_{\max}}
\newcommand{\hdist}[2]{d_\mathbb{H}(#1, #2)}

\newcommand{\petal}{p(\omega)}
\newcommand{\size}[1]{\vert #1 \vert}

\renewcommand{\epsilon}{\varepsilon}

\renewcommand{\geq}{\geqslant}
\renewcommand{\leq}{\leqslant}

\definecolor{lightgray}{rgb}{0.83, 0.83, 0.83}

\title{A new perspective on the Sullivan dictionary \\  via Assouad type dimensions and spectra}
\author{Jonathan M. Fraser and Liam Stuart \\ \\
The University of St Andrews, Scotland\\
E-mails: jmf32@st-andrews.ac.uk and  ls220@st-andrews.ac.uk}

\newlength{\bibitemsep}
\newlength{\bibparskip}
\let\oldthebibliography\thebibliography
\renewcommand\thebibliography[1]{%
  \oldthebibliography{#1}%
  \setlength{\parskip}{\bibitemsep}%
  \setlength{\itemsep}{\bibparskip}%
}

\date{\today}

\begin{document}
\pagenumbering{arabic}
\maketitle
\begin{abstract}
The Sullivan dictionary provides a beautiful correspondence between Kleinian groups acting on hyperbolic space and rational maps of the extended complex plane. An especially direct correspondence exists concerning the dimension theory of the associated limit sets and Julia sets.  In  recent work we established formulae for the Assouad type dimensions and spectra for these fractal sets and certain  conformal measures they support.  This allows a rather more nuanced comparison of the two families in the context of dimension.  In this expository article  we discuss how these results provide new entries in the Sullivan dictionary, as well as revealing striking differences between the two settings. 
 

 \end{abstract}

\section{Introduction}\label{Introduction}
  Seminal work of Sullivan in the 1980s \cite{S2}  resolved a long-standing problem in complex dynamics by proving that the Fatou set of a rational map has no wandering domains. This work served to establish remarkable connections between the dynamics of rational maps and the actions of Kleinian groups.  This connection subsequently stimulated activity in both the complex dynamics and hyperbolic geometry communities and led to what is now known as the \emph{Sullivan dictionary}; see, for example, \cite{MC1}. The Sullivan dictionary provides a framework to study the relationships between Kleinian groups and rational maps.  In many  cases there are analogous results, even with similar proofs, albeit   expressed in a  different language.  See \cite[Table 1]{das} and references therein.

\begin{figure}[H]
\centering
\includegraphics[width=0.9\textwidth]{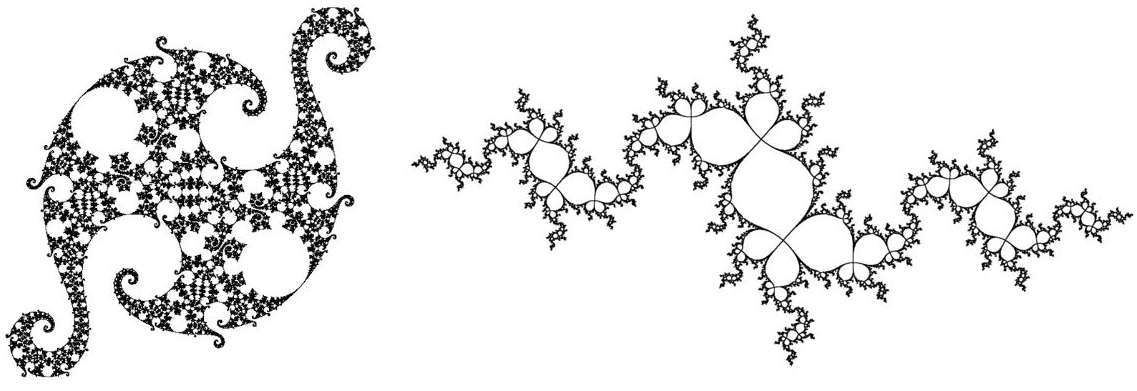}
\caption{Left: an example of a Kleinian limit set. Here $d=2$ and the boundary $\mathbb{S}^2$ has been identified with $\mathbb{R}^2 \cup \{\infty\}$.  Parabolic points with rank 1 are easily identified.  Right: an example of a parabolic Julia set. Parabolic points with petal number 4 are easily spotted. See the following section for definitions and notation. }

\label{juliafig}\label{kleinfig}

\end{figure}

Both  Kleinian groups and rational maps generate important examples of dynamically invariant fractal sets: \emph{limit sets} in the Kleinian case, and \emph{Julia sets} in the rational map case, see Figure \ref{kleinfig}.  The Sullivan dictionary is very well-suited to understanding the connections between these two families of fractal and the correspondence   is especially strong in the context of dimension theory: in both settings there is a `critical exponent' which describes all of the most commonly used notions of fractal dimension (at least in the `geometrically finite' cases).    For Kleinian groups the critical exponent is the Poincar\'e exponent, denoted by $\delta$, and for  rational maps  the critical exponent is   the smallest  zero   of the topological pressure, denoted by $h$. For both geometrically finite Kleinian groups and rational maps the critical exponent coincides with the Hausdorff, packing and box dimensions of the associated fractal as well as the Hausdorff, packing, and entropy dimensions of the associated ergodic conformal measure of maximal dimension.

There has been a recent increase in interest in the \emph{Assouad type dimensions} and these dimensions (and associated dimension spectra) do not behave in such a straightforward manner in the presence of parabolicity.  In particular, the critical exponent does \emph{not necessarily} give the Assouad dimension of the associated fractals. As we shall see, by slightly expanding the family of dimensions considered, a much richer and more varied tapestry of results emerges. In this expository paper we discuss recent work from \cite{Fr1, stuartjulia,stuartkleinian} and show how this can be used to provide a new perspective in the Sullivan dictionary.

\begin{table}[H]
\centering
\begin{tabular}{ c|c } 
Kleinian  & Julia  \\ 
 \hline
 Kleinian group $\Gamma$ & rational map $T$ \\ 
Kleinian limit set $\lset$ &  Julia set $J(T)$ \\
Poincar\'e exponent $\delta$  & critical exponent $h$  \\ 
Patterson-Sullivan measure $\mu$ & $h$-conformal measure $m$ \\
$\haus \lset = \boxd \lset = \delta$ & $\haus J(T) = \boxd J(T) = h$  \\
$\haus \mu =   \delta$ & $\haus m = h$  \\
finite set of inequivalent parabolic points & finite set of parabolic points $\mathbf{\Omega}$ \\
rank of parabolic point $k(p)$ & petal number of parabolic point $\petal$ \\
dimension bound $\delta > k_{\max}/2$ & dimension bound $h > p_{\max}/(1+p_{\max})$ \\

\end{tabular}
\caption{Some well-known `entries' in the Sullivan dictionary.  See the following section for definitions and notation.  In Section \ref{Sull} we describe an expansion of this dictionary, including several new entries as well as some striking differences (`non-entries').}
\end{table}

\section{Definitions and Background}\label{Prelims}
\subsection{Dimensions of sets and measures and `dimension interpolation'}
\label{DimPrelims}
We recall and motivate the key notions from fractal geometry and dimension theory which we use.  For a more in-depth treatment see the books \cite{BP, FK} for background on Hausdorff and box dimensions, and \cite{Fr2} for Assouad type dimensions.    We  work with fractals in two distinct settings.  Kleinian limit sets will be subsets of the $d$-dimensional sphere $\mathbb{S}^d$ which we  view as a subset of $\mathbb{R}^{d+1}$. On the other hand, Julia sets will be subsets of the Riemann sphere  $\hat{\mathbb{C}} = \mathbb{C} \cup \{\infty\}$.  However, by a standard reduction we will assume that the Julia sets are bounded subsets of the complex plane $\mathbb{C}$, which we  identify with $\mathbb{R}^2$.  Therefore, it is convenient to recall dimension theory for bounded subsets of Euclidean space only. 

Let $F \subseteq \mathbb{R}^d$ be non-empty and bounded. Perhaps the most commonly used notion of fractal dimension is the Hausdorff dimension.  We write $\haus F$, $\boxd F$ and $\ubox$  for the \emph{Hausdorff, box} and \emph{upper box dimensions} of $F$, respectively,  but refer the reader to \cite{BP, FK} for the precise definitions.   We write 
\[
\size{F} = \sup_{x,y \in F} |x-y| \in [0,\infty)
\]
 to denote the diameter of  $F$.  Given $r>0$, we write $N_r(F)$ for the smallest number of balls of radius $r$ required to cover $F$. In the last 10 years there has been an increase in interest in the Assouad dimension in the context of fractal geometry. This notion has been of central importance in other fields for much longer, however, and stems from work in embedding theory and conformal geometry, see \cite{MT, R}.   The \textit{Assouad dimension} of $F$ is defined by
\begin{align*}
\aso F = \inf \Bigg\{ s \geq 0 \mid  \exists C>0 \ : \ \forall \ 0<r<R \  :  \ \forall x \in F \  :  \ 
N_r(B(x,R) \cap F) \leq C \left(\frac{R}{r} \right)^{s} \Bigg\}. 
\end{align*}
 The lower dimension is the natural `dual' to the Assouad dimension and it is particularly useful to consider these notions together. The \textit{lower dimension} of $F$ is defined by
\begin{align*}
\low F = \sup \Bigg\{ s \geq 0 \mid  \exists C>0 \ : \ \forall \ 0<r<R \leq |F|  \  :  \ \forall x \in F \  :  \   
N_r(B(x,R) \cap F) \geq C \left(\frac{R}{r} \right)^{s} \Bigg\}
\end{align*}
provided $|F| >0$ and otherwise it is 0. Importantly, for compact $F$ we have 
\[
\low F \leq \haus F  \leq \ubox F \leq \aso F.
\]
 The Assouad and lower spectra were introduced much more recently  in \cite{FYu} and provide an `interpolation' between the box dimension and the Assouad and lower dimensions, respectively.  The motivation for the introduction of these `dimension spectra' was to gain a more nuanced understanding of fractal sets than that provided by the dimensions considered in isolation.  This is already proving a fruitful programme with applications emerging in a variety of  settings including to problems in  harmonic analysis, see work of Anderson, Hughes, Roos and Seeger \cite{AHRS} and \cite{RS}. These spectra provide a parametrised family of dimensions by fixing the relationship between the two scales $r<R$ used to define Assouad and lower dimension.   Studying the dependence on the parameter within this family  thus yields finer and more nuanced information about the local structure of the set.  For example, one may understand which scales `witness' the behaviour described by the Assouad and lower dimensions.  For $\theta \in (0,1)$, the \textit{Assouad spectrum} of $F$ is given by
\begin{align*}
\asospec F = \inf \Bigg\{ s \geq 0 \mid  \exists C>0 \ : \ \forall \ 0<r<1 \  :  \ \forall x \in F \  :  \ 
N_r(B(x,r^{\theta}) \cap F) \leq C \left(\frac{r^{\theta}}{r} \right)^{s} \Bigg\} .
\end{align*}
The \textit{lower spectrum} of $F$, denoted by $\low F$,  is defined similarly by using the parameter to fix the relationship $R=r^\theta$ in the definition of the lower dimension. It was shown in \cite{FYu} that
\begin{align}
\ubox F &\leq \asospec F \leq \min\left\{\aso F, \ \frac{\ubox F}{1-\theta}\right\} \label{basicbound}\\
\low F &\leq \lowspec F \leq \lbox F. \nonumber
\end{align}
In particular, $\asospec F \to \ubox F$ as $\theta \to 0$.  The limit of  $\asospec F$ exists and coincides with the quasi-Assouad dimension. The quasi-Assouad and Assouad dimensions do not necessarily  coincide, but in many cases of interest they do. It is not necessarily true that $\lowspec F \to \lbox F$ as $\theta \to 0$, but it was proved in \cite[Theorem 6.3.1]{Fr2} that this does hold  provided $F$ satisfies a strong form of  dynamical invariance.   Whilst the fractals we study are not quite covered by this result, we shall see that this interpolation holds nevertheless.

There is an analogous dimension theory of measures, and the interplay between the dimension theory of fractals and the measures they support is fundamental to fractal geometry, especially in the dimension theory of dynamical systems.  For example, a problem of interest is to identify dynamical measures witnessing the dimension of the support, e.g. invariant measures of full Hausdorff dimension.  Let $\nu$ be a locally finite Borel measure on $\mathbb{R}^d$, i.e. $\nu(B(x,r)) < \infty$ for all $x \in \mathbb{R}^d$ and $r>0$. We write $\text{supp}(\nu) = \{x \in \mathbb{R}^d \mid \nu(B(x,r)) > 0 \ \text{for all} \ r>0 \}$ for the \textit{support} of $\nu$. We say that $\nu$ is \textit{fully supported} on a set $F \subseteq \mathbb{R}^d$ if $\text{supp}(\nu) = F$. Similar to above, we write $\haus \nu$ for the (lower) Hausdorff dimension of  $\nu$ and note that $\haus \nu \leq \haus \text{supp}(\nu)$ and, for compact $F$, 
\[
\haus F = \sup \{ \haus \nu \mid  \text{supp}(\nu) \subseteq  F\},
\]
 see \cite{MAT}.  The \textit{Assouad dimension} of $\nu$ with $\text{supp}(\nu)=F$ is defined by
\begin{align*}
\aso \nu = \inf \Bigg\{ s \geq 0 \mid  \exists C>0 \ : \ \forall \ 0<r<R< \vert F \vert \  :  \ \forall x \in F \  :  \  \frac{\nu(B(x,R))}{\nu(B(x,r))} \leq C \left(\frac{R}{r} \right)^{s} \Bigg\} 
\end{align*}
and, provided $\vert \text{supp}(\nu) \vert  = |F|> 0$, the \textit{lower dimension} of $\nu$ is given by
\begin{align*}
\low \nu = \sup \Bigg\{ s \geq 0 \mid  \exists C>0 \ : \ \forall \ 0<r<R< \vert F \vert \  :  \ \forall x \in F \  :  \  \frac{\nu(B(x,R))}{\nu(B(x,r))} \geq C \left(\frac{R}{r} \right)^{s} \Bigg\}
\end{align*}
and otherwise it is 0. By convention we assume that $\inf \emptyset = \infty$. The Assouad and lower dimensions of measures  were introduced in \cite{Ka2}, where they were referred to as the upper and lower regularity dimensions, respectively. It is well known (see \cite[Lemma 4.1.2]{Fr2}) that for a Borel probability measure $\nu$ supported on a closed set $F \subseteq \mathbb{R}^{d}$, we have 
\[
\low \nu \leq \low F \leq \aso F \leq \aso \nu
\]
 and, furthermore, we have the  stronger fact that
\[ \aso F = \inf \left\{ \aso \nu \mid \text{$\nu$ is a Borel probability measure fully supported on} \ F\right\}\]
and 
\[ \low F = \sup \left\{ \low \nu \mid \text{$\nu$ is a Borel probability measure fully supported on} \ F\right\}.\]
For $\theta \in (0,1)$, the \textit{Assouad spectrum} of $\nu$, denoted by $\asospec \nu $   and the \textit{lower spectrum} of $\nu$, denoted by $\lowspec \nu $ are defined similarly to the Assouad and lower dimensions but, again, using the parameter $\theta \in (0,1)$ to fix the relationship $R=r^\theta$.  

It is known (see \cite{FFK} for example) that for any measure $\nu$, 
\[
\low \nu \leq \lowspec \nu \leq \asospec \nu \leq \aso \nu
\]
 and, if $\nu$ is fully supported on a closed set $F$, then 
\[
\lowspec \nu \leq \lowspec F \leq \asospec F \leq \asospec \nu.
\]
There are also upper and lower box dimensions for measures, recently introduced in \cite{FFK}.  We omit the formal definitions, referring the reader to \cite{FFK, Fr2}.  Following \cite{FFK}, it is useful to note that 
\[ \ubox F = \inf \left\{ \ubox \nu \mid \text{$\nu$ is a finite Borel measure fully supported on} \ F\right\}\]
with an analogous result for lower box dimension. Furthermore, it was shown that the upper box dimension of $\nu$ can be related to the Assouad spectrum of $\nu$ in a similar manner to sets, that is, for $\theta \in (0,1)$,
\[\ubox \nu \leq \asospec \nu \leq \min\left\{\aso \nu, \frac{\ubox \nu}{1-\theta}\right\}\]
and so $\ubox \nu = \lim_{\theta \to 0} \asospec \nu$.

\subsection{Kleinian groups and limit sets}
\label{KleinPrelims}
For a more thorough study of hyperbolic geometry and Kleinian groups, we refer the reader to \cite{B, M}. For $d \geq 1$, we model $(d+1)$-dimensional hyperbolic space using the Poincar\'e ball model 
\[
\mathbb{D}^{d+1} = \{z \in \mathbb{R}^{d+1} \mid \vert z \vert < 1\}
\]
  equipped with the hyperbolic metric $d_{\mathbb{H}}$ 
and we call the boundary 
\[
\mathbb{S}^d = \{z \in \mathbb{R}^{d+1} \mid \vert z \vert = 1\}
\]
 the \textit{boundary at infinity} of the space $(\mathbb{D}^{d+1},d_{\mathbb{H}})$. We denote by $\text{Con}(d)$ the group of orientation-preserving isometries of  $(\mathbb{D}^{d+1},d_{\mathbb{H}})$. We say that a group is \textit{Kleinian} if it is a discrete subgroup of $\text{Con}(d)$,
and given a Kleinian group $\Gamma$, the \textit{limit set} of $\Gamma$ is defined to be  $\lset = \overline{\Gamma(\mathbf{0})} \setminus \Gamma(\mathbf{0})$ where $\mathbf{0} = (0,\dots,0) \in \mathbb{D}^{d+1}$. It is well known that $\lset$ is a compact $\Gamma$-invariant subset of $\mathbb{S}^d$, see Figure \ref{kleinfig}. If $\lset$ contains zero, one or two points, it is said to be \textit{elementary}, and otherwise it is \textit{non-elementary}. In the non-elementary  case, $\lset$ is a perfect set, and often has a complicated fractal structure. We consider \emph{geometrically finite} Kleinian groups.  Roughly speaking, this means that there is a fundamental domain with finitely many sides but we refer the reader to \cite{BO} for a precise definition. We define the \textit{Poincar\'e exponent} of a Kleinian group $\Gamma$ to be 
\[\delta = \inf\left\{s>0 \mid \sum_{g \in \Gamma} e^{-s \hdist{\mathbf{0}}{g(\mathbf{0})}} < \infty \right\}.\]
Due to work of  Patterson and Sullivan \cite{Pa1, S1}, it is known that for a non-elementary geometrically finite Kleinian group $\Gamma$, the Hausdorff dimension of the limit set is equal to $\delta$. It was proved independently by Bishop and Jones \cite[Corollary 1.5]{BJ} and Stratmann and Urba\'nski \cite[Theorem 3]{SU3} that the box and packing dimensions of the limit set are also equal to $\delta$.  Even in the non-elementary geometrically \emph{infinite} case, $\delta$ is still an important quantity.  In fact it  always gives the Hausdorff dimension of the \emph{radial} limit set, and therefore also provides a lower bound for the Hausdorff dimension of the limit set, see \cite{BJ}. 

From now on we only discuss the non-elementary geometrically finite case. We write $\mu$ to denote the \textit{Patterson-Sullivan measure}, which is a measure first constructed by Patterson in \cite{Pa1}.  Strictly speaking there  is a family of (mutually equivalent) Patterson-Sullivan measures.  However, we may fix one for simplicity (and hence talk about \emph{the} Patterson-Sullivan measure since the dimension theory is the same for each measure).  The geometry of $\Gamma$, $\lset$ and $\mu$ are heavily related.  For example,  $\mu$ is a conformal $\Gamma$-ergodic Borel probability measure which is fully  supported on $\lset$.  Moreover, $\mu$ has Hausdorff, packing and entropy dimension equal to $\delta$, see \cite{SV}.  The limit set is $\Gamma$-invariant in the strong sense that $g(\lset) = \lset$ for all $g \in \Gamma$.  However, $\mu$ is only quasi-invariant and $\mu \circ g$ is related to $\mu$ by a geometric transition rule, see \cite[Chapter 14]{borth} for a more detailed exposition of this. 

If $\Gamma$ contains no parabolic elements, then \[\aso \lset = \low \lset = \aso \mu = \low \mu = \boxd \mu= \delta,\] 
 see \cite{Fr1}.  Therefore, we  assume  from now on that $\Gamma$ contains at least one parabolic element.

Let  $P \subseteq \lset$ denote the countable set of parabolic fixed  points. For  $p \in P$ write $k(p)$ to denote the maximal rank of a free abelian subgroup of the stabiliser of $p$ (in $\Gamma$) and call this the \textit{rank} of $p$. We write 
\begin{align*}
\kmin &= \min\{k(p) \mid p \in P\}\\
\kmax &= \max\{k(p) \mid p \in P\}.
\end{align*}
 It was proven in \cite{S1} that $\delta > \kmax/2$.

\subsection{Rational maps and Julia sets}
\label{JuliaPrelims}
For a more detailed discussion of the dynamics of rational maps, see \cite{CG, milnor}. Let $T:\hat{\mathbb{C}} \rightarrow \hat{\mathbb{C}}$ denote a rational map of degree at least 2, and write  $J(T)$ to denote  the \textit{Julia set} of $T$, which is equal to the closure of the repelling periodic points of $T$, see Figure \ref{juliafig}.  The Julia set is closed and $T$-invariant.  We may assume that $J(T)$ is a compact subset of $\mathbb{C}$ by a standard reduction. 

 A periodic point $\xi \in  \mathbb{C}$ with period $p$ is said to be \textit{rationally indifferent} (or \emph{parabolic}) if
$(T^p)^{'}(\xi) = e^{2 \pi i q}$ for some $q \in \mathbb{Q}$. We say that $T$ and $J(T)$ are \textit{parabolic} if $J(T)$ contains no critical points of $T$, but contains at least one parabolic point.  Define $h$ to be the  smallest zero of the topological pressure  $t \mapsto P(T, -t\text{log} \size{T'})$.   In the parabolic setting, it was proven in \cite{DU2} that     $\haus J(T) = h$.  Furthermore, in \cite{DU4} it was shown that the box and packing dimensions of $J(T)$ are equal to $h$.  Due to work of Aaronson, Denker and Urba\'nski \cite{ADU, DU1, DU2} it is known that, for parabolic $T$, there exists a unique atomless $h$-conformal probability measure $m$ supported on  $J(T)$.  It also follows from, for example, \cite{SU2} that $m$ is exact dimensional and therefore the Hausdorff,  packing and entropy dimensions of $m$ are also given by $h$.  

 If $T$ contains no critical points nor parabolic points, then it is \emph{hyperbolic} and, analogous to case of  geometrically finite  Kleinian groups with no parabolic elements,
\[\aso J(T) = \low J(T) = \aso m = \low m =\boxd m = h,\]
see \cite{stuartjulia}.  Therefore, we assume from now on that $T$ is parabolic.  

Write  $\mathbf{\Omega}$ to denote the finite set of parabolic points of $T$, and let
\[\mathbf{\Omega}_0 = \{\xi \in \mathbf{\Omega} \mid T(\xi) = \xi, \ T'(\xi) = 1\}.\]
As $J(T^{n}) = J(T)$ for every $n \in \mathbb{N}$, we may assume without loss of generality that $\mathbf{\Omega} = \mathbf{\Omega_0}.$   Following \cite{DU4, SU1},   for each $\omega \in \mathbf{\Omega}$, we can find a ball $U_{\omega} = B(\omega,r_{\omega})$ with sufficiently small radius such that on $B(\omega,r_{\omega})$, there exists a unique holomorphic inverse branch $T_{\omega}^{-1}$ of $T$ such that $T_{\omega}^{-1}(\omega) = \omega$. For a parabolic point $\omega \in \mathbf{\Omega}$, the Taylor series of $T$ about $\omega$ is of the form
\[z + a(z-\omega)^{\petal + 1} + \cdots.\]
We call $\petal$ the \textit{petal number} of $\omega$, and we write
\begin{align*}
\pmin &= \min\{p(\omega) \mid \omega \in \mathbf{\Omega}\}\\
\pmax &= \max\{p(\omega) \mid \omega \in \mathbf{\Omega}\}.
\end{align*} 
It was proven in \cite{ADU} that $h > \pmax/(1+\pmax)$.

\section{A new perspective on the Sullivan dictionary} \label{Sull}

\subsection{Recent results on Assouad type dimensions and spectra}\label{Results}\label{KleinResults} \label{JuliaResults}

In this subsection we state various recent results concerning geometrically finite Kleinian groups and parabolic Julia sets which provide a new perspective on the Sullivan dictionary in the context of dimension theory.  The Assouad and lower dimensions of limit sets of geometrically finite Kleinian groups and associated Patterson-Sullivan measures were found in \cite{Fr1}.  The analogous results for Julia sets were proved in \cite{stuartjulia}.  The results concerning Assouad type spectra  were proved in \cite{stuartjulia, stuartkleinian}. Throughout we fix $\theta \in (0,1)$.

\subsubsection{Patterson-Sullivan measure $\mu$}

\begin{align*}
\normalfont{\aso} \mu &= \max\{2\delta-k_{\min},k_{\max}\}\\
\normalfont{\boxd} \mu &= \max\{2\delta-\kmin, \delta\} \\
\normalfont{\low} \mu &= \min\{2\delta-k_{\max},k_{\min}\}
\end{align*}
\begin{align*}
\normalfont{\asospec} \mu &=  \Bigg\{ 
    \begin{array}{cc}
       \delta + \min\left\{1,\frac{\theta}{1-\theta}\right\}(\kmax-\delta) & \delta < \kmin \\
2\delta -\kmin + \min\left\{1,\frac{\theta}{1-\theta}\right\}(\kmin+\kmax-2\delta)  & \kmin \leq \delta < (\kmin+\kmax)/2 \\
       2\delta-\kmin & \delta \geq (\kmin+\kmax)/2
    \end{array} \\
\normalfont{\lowspec} \mu &= \Bigg\{ 
    \begin{array}{cc}
             2\delta-\kmax & \delta \leq (\kmin+\kmax)/2 \\
2\delta-\kmax - \min\left\{1,\frac{\theta}{1-\theta}\right\}(2\delta-\kmin-\kmax)  & (\kmin+\kmax)/2 < \delta \leq \kmax \\
 \delta - \min\left\{1,\frac{\theta}{1-\theta}\right\}(\delta-\kmin) & \delta > \kmax  
    \end{array}
\end{align*}

 \subsubsection{Kleinian limit sets $\lset$}
\begin{align*}
\normalfont{\aso} \lset &= \max\{\delta,k_{\max}\} \\
\normalfont{\low} \lset  &= \min\{\delta,k_{\min}\}
\end{align*}
\begin{align*}
\normalfont{\asospec} \lset  &=  \Bigg\{ 
    \begin{array}{cc}
     \delta + \min\left\{1,\frac{\theta}{1-\theta}\right\}(\kmax-\delta) & \delta < \kmax\\
\delta  & \delta \geq \kmax
    \end{array}\\
\normalfont{\lowspec} \lset  &=  \Bigg\{ 
    \begin{array}{cc}
     \delta   &\delta \leq \kmin\\
\delta - \min\left\{1,\frac{\theta}{1-\theta}\right\}(\delta-\kmin)  & \delta > \kmin
    \end{array}
\end{align*}

 \subsubsection{$h$-conformal measures  $m$}

\begin{align*}
\normalfont{\aso} m &= \max\{1,h+(h-1)\pmax\}\\
\normalfont{\boxd} m &=  \max\{h,h+(h-1)\pmax\}\\
\normalfont{\low} m &= \min\{1,h+(h-1)\pmax\}
\end{align*}
\begin{align*}
\normalfont{\asospec} m  &=  \Bigg\{ 
    \begin{array}{cc}
     h+\min\left\{1,\frac{\theta \, \pmax}{1-\theta}\right\}(1-h)  & h<1\\
h+(h-1)\pmax & h\geq 1
    \end{array}\\
\normalfont{\lowspec} m  &=  \Bigg\{ 
    \begin{array}{cc}
    h+(h-1)\pmax  & h<1\\
h+\min\left\{1,\frac{\theta \, \pmax}{1-\theta}\right\}(1-h) & h\geq 1
    \end{array}
\end{align*}

 \subsubsection{Julia sets $J(T)$}

\begin{align*}
\normalfont{\aso} J(T) &= \max\{1,h\}\\
\normalfont{\low} J(T) &= \min\{1,h\}
\end{align*}
\begin{align*}
\normalfont{\asospec} J(T)  &=  \Bigg\{ 
    \begin{array}{cc}
  h+\min\left\{1,\frac{\theta \, \pmax}{1-\theta}\right\}(1-h)  & h<1\\
h  & h\geq 1
    \end{array}\\
\normalfont{\lowspec} J(T)  &=  \Bigg\{ 
    \begin{array}{cc}
  h   & h<1\\
 h+\min\left\{1,\frac{\theta \, \pmax}{1-\theta}\right\}(1-h)  & h\geq 1
    \end{array}
\end{align*}

\subsection{New entries in the Sullivan dictionary}

Given the array of results in the previous section, it is clear that there are some parallels between the Kleinian and Julia settings akin to the Sullivan dictionary.  Here we take a closer look at some of these parallels.

 1) \textit{Interpolation between dimensions.}  In both settings, the Assouad spectrum always interpolates between the upper box and Assouad dimensions of the respective sets and measures regardless of what form it takes, that is,  $\displaystyle \lim_{\theta \to 1} \asospec F = \aso F$ where $F$ can be replaced by $\mu, \lset, m$ or $J(T)$. Recall that this interpolation does  not  hold in general. Similar interpolation holds as $\theta \to 1$ for the lower dimensions and spectra.   

2) \textit{Failure to witness the box dimension of measures.}  For the measures $\mu$ and $m$, the lower spectrum does not generally tend to the box dimension as $\theta \to 0$.  In fact, if the lower spectrum does tend to the box dimension as $\theta \to 0$, then it is constant and  $\delta = \kmin = \kmax$ (in the Kleinian setting) and $h=1$ (in the Julia setting).

 3) \textit{General form of the spectra.}    For  $F$ a given set or measure, consider 
\[\rho = \inf\{\theta \in (0,1) \mid \asospec F = \aso F\}.\]
Following some algebraic manipulation, we find that in all cases
\begin{equation} \label{3pf}
\asospec F = \min\left\{ \boxd F + \frac{(1-\rho)\theta}{(1-\theta)\rho} (\aso F - \boxd F), \aso F \right\} 
\end{equation}
where $F$ can be replaced by $\mu, \lset, m$ or $J(T)$.  This formula, and the fact that the    Assouad spectrum can be expressed purely  in terms of the `phase transition' $\rho$ together with the box  and Assouad dimensions, has appeared in a variety of settings, see \cite[Section 17.7]{Fr2} and the discussion therein.  For example, this formula also holds for self-affine Bedford-McMullen carpets.  The phase transition $\rho$ often has a natural `geometric significance' for the objects involved and opens the door to a new `dictionary' extending beyond the setting discussed here.  It is worth noting that \eqref{3pf} does not hold generally, even failing for simple  examples such as  the elliptical spirals considered in \cite{burrell}. 

 4) \textit{The phase transition and the Hausdorff dimension bound.}  There is a correspondence between the phase transition $\rho$  and the general lower bounds for the Hausdorff dimension.   Applying   \eqref{basicbound} shows that, for any non-empty bounded set $F$,  $\rho \geq 1 - \ubox F/\aso F.$  When the spectra are non-constant,  in the Kleinian setting we always have $\rho = 1/2$, and in the Julia setting we always have $\rho = 1/(1+\pmax)$.  Combining this with  the general Hausdorff dimension bounds  $\delta > \kmax/2 = \kmax \rho$ and $h > \pmax/(1+\pmax) = \pmax \rho$ in both settings yields $\rho > 1 - \ubox F/\aso F$, showing that the upper bound from \eqref{basicbound} is never achieved in either  setting (but is asymptotically sharp).  

 5) \textit{The realisation problem.} Given the interplay between dimensions of sets and dimensions of measures seen in Section \ref{DimPrelims}, one may ask if it is possible to construct an (invariant, or quasi-invariant) measure $\nu$ which \textit{realises} the dimensions of an (invariant) set $F$, that is, $\dim \nu = \dim F$. One can ask this about a particular choice of dimension $\dim$ or if a single measure can be constructed to solve the problem for several notions of dimension simultaneously. We note that the measures $\mu$ and $m$ \emph{always} realise the Hausdorff dimensions of $\lset$ and $J(T)$ respectively. As for the Assouad and lower dimensions,   $\mu$ realises the Assouad dimension of $\lset$ when $\delta \leq (\kmin+\kmax)/2$ and realises the lower dimension when $\delta \geq (\kmin+\kmax)/2$. Similarly, for $m$ to realise the Assouad dimension of $J(T)$ we require $h \leq 1$, and for $m$ to realise the lower dimension of $J(T)$ we require $h \geq 1$. A similar relationship holds for the box dimension too: in the Kleinian setting we require $\delta \leq \kmin$ and in the Julia setting we require $h \leq 1$.

 6) \textit{A  special case.} Finally, we observe that in the (very) special case $\kmin = \kmax = \pmax = 1$, the formulae for the Assouad type dimensions and spectra are identical in the Kleinian and Julia settings.  Does this suggest that this special case is one where we can expect the Sullivan dictionary to yield a particularly strong correspondence in other settings?

\subsection{New non-entries in the Sullivan dictionary}

Here we discuss some notable differences between the Kleinian and Julia settings.  These are especially interesting to us since the Sullivan dictionary previously provided a very strong parallel in the context of dimension theory.

 1) \textit{Assouad dimension.} Our results show that Julia sets of parabolic rational maps can never have full Assouad dimension, that is, we always have $\aso J(T) <2$. This uses our result together with \cite[Theorem 8.8]{ADU} which proves that $h<2$. This is in stark contrast to the situation for Kleinian limit sets where it is perfectly possible for the Assouad dimension to be full, that is, $\Gamma \in \text{Con}(d)$ with $\aso \lset = d = \dim \mathbb{S}^d$ for any integer $d \geq 1$.   This can even happen when the limit set is nowhere dense (that is, when $\delta<d$, see  \cite[Theorem D]{Tuk}). We note that   $\aso J(T) <2$ also follows from \cite[Theorem 1.4]{LG}, where it was proved that  parabolic Julia sets are porous, together with  \cite[Theorem 5.2]{LK}, which shows  that porous sets in $\mathbb{R}^d$ must have Assouad dimension strictly less than $d$. Our results can thus be viewed as a refinement of the observation that parabolic Julia sets are porous.  We note that Julia sets of general  rational maps need not be porous, and may even have positive area.  This was proved to be possible even within the  quadratic family by  Buff and Ch\'eritat \cite{buff}. We  proved in \cite{stuartjulia} that Julia sets with Cremer fixed points have Assouad dimension 2 (and are therefore not porous).

 2) \textit{Lower dimension.} Our results, together with the standard  bound  $h > \pmax/(1+\pmax)$,  show that $\low J(T) = \min\{1,h\} > \pmax/(1+\pmax)$, that is, the lower dimension respects the general lower bound satisfied by the Hausdorff dimension.  Again, this is in stark contrast to the situation for Kleinian limit sets where the standard bound for Hausdorff dimension is $\delta > \kmax/2$ but  $\low \lset = \min\{\kmin, \delta \} \leq  \kmax/2$ is possible, even  in the $d=2$ case.

 3) \textit{Relationships between dimensions.} An interesting aspect of dimension theory is to consider what configurations are possible between the different notions of dimension in a particular setting.  We refer the reader to \cite[Section 17.5]{Fr2} for a more general discussion of this.  Our results show that
\[\low J(T) < \haus J(T) < \aso J(T)\]
is impossible in the Julia setting but the analogous configuration  \emph{is} possible in the Kleinian setting, even  in the $d=2$ case.
 
\begin{table}[H]
\centering
\begin{tabular}{ c|c|c|c } 
Configuration & Fuchsian & Kleinian & Julia  \\ 
 \hline
 L=H=A  & $\checkmark$ & $\checkmark$ & $\checkmark$\\ 
L=H$<$A  & $\checkmark$  &  $\checkmark$ & $\checkmark$\\ 
L$<$H=A & $\times$ & $\checkmark$ & $\checkmark$ \\
L$<$H$<$A & $\times$ &   $\checkmark$ & $\times$
\end{tabular}
\caption{Summarising the possible relationships between the lower, Hausdorff,  and Assouad dimensions of geometrically finite Fuchsian limit sets, geometrically finite Kleinian limits sets and parabolic Julia sets with the obvious labelling. The label `Fuchsian' refers to the Kleinian setting when $d=1$. The symbol $\checkmark$ means that the configuration is possible, and $\times$ means the configuration is  impossible.  In other situations it is interesting to add box dimension into this discussion, but here this always coincides with Hausdorff dimension and so we omit it.}
\end{table}

4) \textit{Form of the spectra.} Turning our attention to  measures, the Assouad and lower spectra of $\mu$ in the Kleinian setting can take 3 different forms, in comparison to the Julia setting where we only have 2 possibilities for $m$. Furthermore, in the Kleinian setting, both $\kmin$ and $\kmax$ appear in the formulae for the Assouad and lower spectra, sometimes simultaneously, but in the Julia setting only $\pmax$ appears.

5) \textit{The realisation problem for dimension spectra.} One can also extend the realisation problem to the Assouad and lower spectra: when does an (invariant) set support an (invariant, or quasi-invariant) measure with equal Assouad or lower spectra?  In the Kleinian setting, we have $\asospec \mu = \asospec \lset$ when $\delta \leq \kmin$ and $\lowspec \mu = \lowspec \lset$ when $\delta \geq \kmax$. This can leave a gap when $\kmin < \delta < \kmax$ where neither of the spectra are realised by the Patterson-Sullivan measure. 
This is in contrast to the Julia setting where $\asospec m = \asospec J(T)$ when $h \leq 1$ and $\lowspec m = \lowspec J(T)$  when $h \geq 1$, and so at least one of the spectra is always realised by $m$.

6) \textit{Dimension spectra as a fingerprint.} Suppose it is \emph{not} true that $\kmin = \kmax   = \pmax = 1$.  Then simply by looking at plots of the Assouad and lower spectra, one can determine whether the set in question is a Kleinian limit set or a Julia set.    Whenever the Assouad   spectrum is non-constant in either the Kleinian or Julia setting, there is a unique phase transition at 
\[
\rho = \inf\{\theta \in (0,1) \mid \asospec F = \aso F\}.
\]
  However,  $\rho=1/2$ in the Kleinian setting and $\rho = 1/(1+\pmax)$ in the Julia setting.    Note that in the Kleinian setting the phase transition is  constant across all Kleinian limit sets, whereas in the Julia setting the phase transition depends on the rational map $T$.    This allows one to distinguish between the Assouad spectrum of a Kleinian limit set and a Julia set just by looking at the phase transition, provided $\pmax \neq 1$.  However, even if $\pmax = 1$, the spectra will still distinguish between the two settings provided we do not also have  $\kmin = \kmax   = 1$.  

\subsection{Examples}\label{Examples}
We plot the Assouad and lower spectra for some examples. In the Kleinian setting, we assume that $d=2$ throughout for a more direct comparison with the Julia setting, and plot the following cases:  $\delta < \kmin$,  $\delta > \kmax$, and $\kmin < \delta < \kmax$. In the Julia setting, we plot examples with: $h < 1$ and  $h > 1$. The following are plots of  the Assouad and lower spectra as functions of $\theta \in (0,1)$. The spectra of $\mu$ and $m$ are plotted with dashed lines, and the spectra of $\lset$ and $J(T)$ by solid lines. The Assouad spectra are plotted in black and the lower spectra are plotted in grey.  
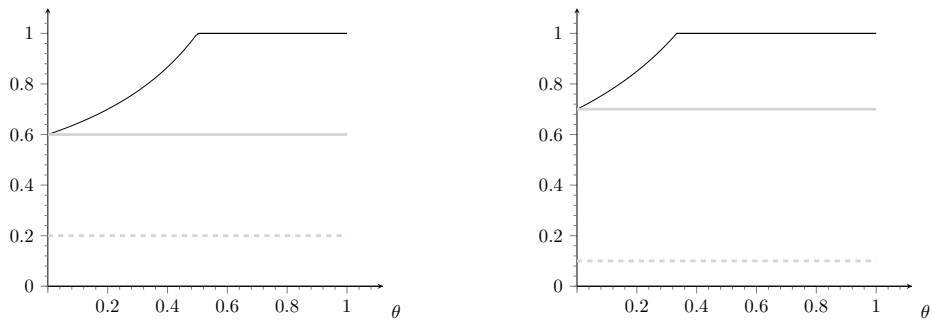
\begin{figure}[H]
\centering

\begin{tikzpicture}[scale=0.65]
\begin{axis}[
    axis lines = left,
    tick align = outside,
    xlabel = $\theta$,
    x label style={at={(axis description cs:1.04,0.05)}},
    y label style={at={(axis description cs:0.17,1.12)},anchor=north,rotate=270},
    extra x ticks={0.00,0.04,0.08,0.12,0.16,0.24,0.28,0.32,0.36,0.44,0.48,0.52,0.56,0.64,0.68,0.72,0.76,0.84,0.88,0.92,0.96,1.04,1.08},
    extra y ticks={0.0,0.04,0.08,0.12,0.16,0.24,0.28,0.32,0.36,0.44,0.48,0.52,0.56,0.64,0.68,0.72,0.76,0.84,0.88,0.92,0.96,1.04},
    extra x tick labels={},
    extra y tick labels={},
    extra x tick style={major tick length=2pt},
    extra y tick style={major tick length=2pt},
]

\addplot [
    domain=0:1, 
    samples=100, 
    color=black,
]
{0.6+0.4*min(1,(x/(1-x)))};

\addplot [line width=1.6pt,
    domain=0:1, 
    samples=100, 
    color=lightgray ,
]
{0.6};

\addplot [line width=1.6pt,
    domain=0:1, 
    samples=100, 
    color=lightgray ,
    style=dashed
]
{0.2};

\addplot [
    domain=0:1.12, 
    samples=100, 
    color=white,
]
{1.1};

\addplot [
    domain=0:1.12, 
    samples=100, 
    color=black,
]
{0};

\addplot [
    domain=0.5:1, 
    samples=100, 
    color=black,
]
{1};
\end{axis}
\draw[fill,white] (0,-0.35) circle (0.2cm);
\end{tikzpicture} \qquad \qquad 
\begin{tikzpicture}[scale=0.65]
\begin{axis}[
    axis lines = left,
    tick align = outside,
    xlabel = $\theta$,
    x label style={at={(axis description cs:1.04,0.05)}},
    y label style={at={(axis description cs:0.17,1.12)},anchor=north,rotate=270},
    extra x ticks={0.00,0.04,0.08,0.12,0.16,0.24,0.28,0.32,0.36,0.44,0.48,0.52,0.56,0.64,0.68,0.72,0.76,0.84,0.88,0.92,0.96,1.04,1.08},
    extra y ticks={0.0,0.04,0.08,0.12,0.16,0.24,0.28,0.32,0.36,0.44,0.48,0.52,0.56,0.64,0.68,0.72,0.76,0.84,0.88,0.92,0.96,1.04},
    extra x tick labels={},
    extra y tick labels={},
    extra x tick style={major tick length=2pt},
    extra y tick style={major tick length=2pt},
]

\addplot [
    domain=0:1, 
    samples=100, 
    color=black,
]
{0.7+0.3*min(1,((2*x)/(1-x)))};

\addplot [line width=1.6pt,
    domain=0:1, 
    samples=100, 
    color=lightgray,
]
{0.7};

\addplot [line width=1.6pt,
    domain=0:1, 
    samples=100, 
    color=lightgray,
    style=dashed,
]
{0.1};

\addplot [
    domain=0:1.12, 
    samples=100, 
    color=white,
]
{1.1};

\addplot [
    domain=0:1.12, 
    samples=100, 
    color=black,
]
{0};

\addplot [
    domain=0.333:1, 
    samples=100, 
    color=black,
]
{1};
\end{axis}
\draw[fill,white] (0,-0.35) circle (0.2cm);
\end{tikzpicture}
\caption{Left: a Kleinian limit set with $\delta=0.6$ and $\kmin=\kmax=1$. Right:  a Julia set with $h=0.7$ and $\pmax=2$.}
\end{figure}
\vspace{-0.4cm}

\begin{figure}[H]
\centering

\begin{tikzpicture}[scale=0.65]
\begin{axis}[
    axis lines = left,
    tick align = outside,
    xlabel = $\theta$,
    x label style={at={(axis description cs:1.04,0.05)}},
    y label style={at={(axis description cs:0.17,1.12)},anchor=north,rotate=270},
    extra x ticks={0.00,0.04,0.08,0.12,0.16,0.24,0.28,0.32,0.36,0.44,0.48,0.52,0.56,0.64,0.68,0.72,0.76,0.84,0.88,0.92,0.96,1.04,1.08},
    extra y ticks={0.0,0.2,0.4,0.6,0.8,1.2,1.4,1.6,1.8,2.2,2.4,2.6,2.8,3.2},
    extra x tick labels={},
    extra y tick labels={},
    extra x tick style={major tick length=2pt},
    extra y tick style={major tick length=2pt},
]

\addplot [
    domain=0:1, 
    samples=100, 
    color=black,
]
{1.9};

\addplot [
    domain=0:1, 
    samples=100, 
    color=black,
    style=dashed,
]
{2.8};

\addplot [line width=1.6pt,
    domain=0:1, 
    samples=100, 
    color=lightgray ,
]
{1.9-0.9*min(1,(x/(1-x)))};

\addplot [
    domain=0:1.12, 
    samples=100, 
    color=white,
]
{3.4};

\addplot [
    domain=0:1.12, 
    samples=100, 
    color=black,
]
{0};

\addplot [
    domain=0.5:1, 
    samples=100, 
    color=lightgray,
]
{1};
\end{axis}
\draw[fill,white] (0,-0.35) circle (0.2cm);
\end{tikzpicture} \qquad \qquad 
\begin{tikzpicture}[scale=0.65]
\begin{axis}[
    axis lines = left,
    tick align = outside,
    xlabel = $\theta$,
    x label style={at={(axis description cs:1.04,0.05)}},
    y label style={at={(axis description cs:0.17,1.12)},anchor=north,rotate=270},
    extra x ticks={0.00,0.04,0.08,0.12,0.16,0.24,0.28,0.32,0.36,0.44,0.48,0.52,0.56,0.64,0.68,0.72,0.76,0.84,0.88,0.92,0.96,1.04,1.08},
    extra y ticks={0.0,0.2,0.4,0.6,0.8,1.2,1.4,1.6,1.8,2.2,2.4,2.6,2.8,3.2},
    extra x tick labels={},
    extra y tick labels={},
    extra x tick style={major tick length=2pt},
    extra y tick style={major tick length=2pt},
]

\addplot [
    domain=0:1, 
    samples=100, 
    color=black,
]
{1.4};

\addplot [
    domain=0:1, 
    samples=100, 
    color=black,
    style=dashed,
]
{3};

\addplot [line width=1.6pt,
    domain=0:1, 
    samples=100, 
    color=lightgray,
]
{1.4-0.4*min(1,((4*x)/(1-x)))};

\addplot [
    domain=0:1.12, 
    samples=100, 
    color=white,
]
{3.4};

\addplot [
    domain=0:1.12, 
    samples=100, 
    color=black,
]
{0};

\addplot [
    domain=0.2:1, 
    samples=100, 
    color=lightgray,
]
{1};
\end{axis}
\draw[fill,white] (0,-0.35) circle (0.2cm);
\end{tikzpicture}
\caption{Left: a Kleinian limit set with $\delta=1.9$ and $\kmin=\kmax=1$. Right:  a Julia set with $h=1.4$ and $\pmax=4$.}
\end{figure}
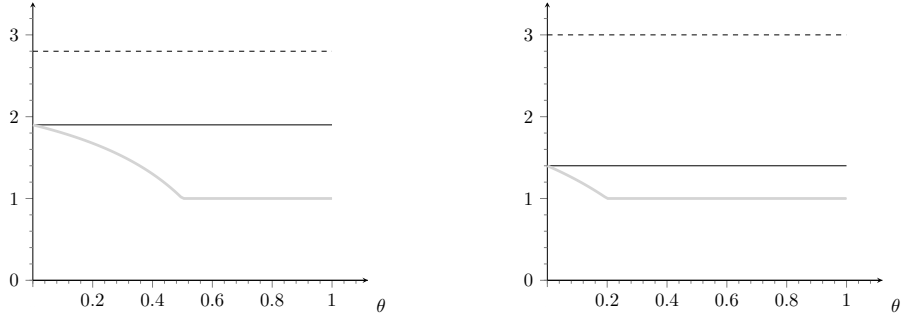

\vspace{-0.4cm}

\begin{figure}[H]
\centering
\begin{tikzpicture}[scale=0.65]
\begin{axis}[
    axis lines = left,
    tick align = outside,
    xlabel = $\theta$,
    x label style={at={(axis description cs:1.04,0.05)}},
    y label style={at={(axis description cs:0.17,1.12)},anchor=north,rotate=270},
    extra x ticks={0.00,0.04,0.08,0.12,0.16,0.24,0.28,0.32,0.36,0.44,0.48,0.52,0.56,0.64,0.68,0.72,0.76,0.84,0.88,0.92,0.96,1.04,1.08},
    extra y ticks={0.0,0.2,0.4,0.6,0.8,1.2,1.4,1.6,1.8,2.2,2.4,2.6,2.8,3.2},
    extra x tick labels={},
    extra y tick labels={},
    extra x tick style={major tick length=2pt},
    extra y tick style={major tick length=2pt},
]

\addplot [
    domain=0:1, 
    samples=100, 
    color=black,
]
{1.7+0.3*min(1,(x/(1-x)))};

\addplot [
    domain=0:1, 
    samples=100, 
    color=black,
    style=dashed,
]
{2.4};

\addplot [line width=1.6pt,
    domain=0:1, 
    samples=100, 
    color=lightgray ,
]
{1.7-0.7*min(1,(x/(1-x)))};

\addplot [line width=1.6pt,
    domain=0:1, 
    samples=100, 
    color=lightgray ,
    style=dashed,
]
{1.4-0.4*min(1,(x/(1-x)))};

\addplot [
    domain=0:1.12, 
    samples=100, 
    color=white,
]
{3.4};

\addplot [
    domain=0:1.12, 
    samples=100, 
    color=black,
]
{0};

\addplot [
    domain=0.5:1, 
    samples=100, 
    color=black,
]
{2};

\addplot [
    domain=0.5:1, 
    samples=100, 
    color=lightgray,
]
{1};
\end{axis}
\draw[fill,white] (0,-0.35) circle (0.2cm);
\end{tikzpicture}
\caption{A Kleinian limit set with $\delta = 1.7$, $\kmin=1$ and $\kmax=2$. In the Julia setting we always have either $\asospec m = \asospec J(T)$ or $\lowspec m = \lowspec J(T)$, and so plots of this form are impossible in the Julia setting.}
\end{figure}
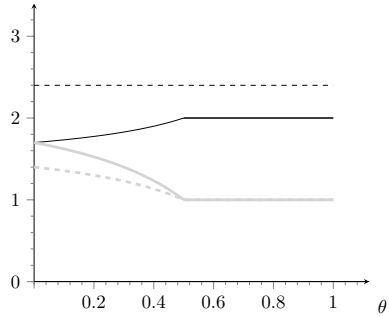

\begin{center} \textbf{Acknowledgements}
\end{center}
\vspace{-3mm}
 JMF was financially supported by an \textit{EPSRC Standard Grant} (EP/R015104/1) and a \textit{Leverhulme Trust Research Project Grant} (RPG-2019-034). LS was financially supported by the University of St Andrews.

\bibliographystyle{apalike}
\addcontentsline{toc}{section}{References}
\bibliography{Sullivan}
\end{document}